\documentclass[journal]{IEEEtran}

\usepackage[hyphens]{url}
\usepackage{hyperref}
\usepackage{breakurl}
\usepackage[T1]{fontenc}

\DeclareTextSymbol{\degre}{T1}{6}
\DeclareTextSymbol{\degre}{OT1}{23}
\usepackage[cmex10]{amsmath}
\usepackage{amsfonts}

\usepackage{dsfont}
\usepackage{graphicx}
\usepackage[noadjust]{cite}
\usepackage{bm}
\usepackage[usenames, dvipsnames]{xcolor} 
\usepackage{tikz}

\pgfdeclareimage{ImageEV}{EV2.jpg}
\usetikzlibrary{calc}
\usetikzlibrary{plotmarks}
\usetikzlibrary{shapes,arrows}
\usepackage[ruled]{algorithm2e}
\usepackage{graphics}
\usepackage{textcomp}
\usepackage{amssymb}
\usepackage{amscd}
\usepackage{epsfig}
\usepackage{psfrag}
\usepackage{rotating}
\usepackage{color}
\usepackage{eurosym}

\usepackage{amsthm}
\usepackage{graphics,color,booktabs,multirow}
\usepackage{amsmath,graphicx,graphics,amsfonts,color}
\usepackage{setspace}
\usepackage{amssymb}
\usepackage{bm}
\usepackage{subfig}
\usepackage{inputenc}
\usepackage{tikz}
\usepackage{bm}
\usepackage{etex}
\usepackage{pstricks}
\usepackage{pstricks-add}
\usepackage{pst-plot}
\usepackage{pst-math}
\usepackage{dsfont}
\usepackage{slashbox}
\newcommand{\mb}{\mathbf}
\newcommand{\mr}{\mathrm}

\newcommand{\bs}{\boldsymbol}
\newcommand{\mc}{\mathcal}

\newcommand{\tcb}{\textcolor{black}}

\newcommand{\ds}{\displaystyle}

\newtheorem{theorem}{Theorem}[section]
\newtheorem{definition}{Definition}

\newtheorem{proposition}[theorem]{Proposition}

\begin{document}
\title{Reducing the Impact of EV Charging Operations on the Distribution Network}
\author{Olivier Beaude, Samson Lasaulce, Martin Hennebel, and Ibrahim Mohand-Kaci
\thanks{O. Beaude is with Renault SAS, L2S and GeePs - CentraleSup\'elec; S. Lasaulce is with CNRS; M. Hennebel is with GeePs - CentraleSup\'elec and I. Mohand-Kaci is with Renault SAS; email: $\lbrace$beaude, lasaulce$\rbrace$@lss.supelec.fr, martin.hennebel@centralesupelec.fr, ibrahim.mohand-kaci@renault.com. A fraction of the material in this paper was presented in \cite{Beaude2012}.}}
\markboth{Transactions on Smart Grid}%
{Shell \MakeLowercase{\textit{et al.}}: Bare Demo of IEEEtran.cls for Journals}

\maketitle

\begin{abstract}
A key assumption made in this paper is that electric vehicle (EV) battery charging profiles are rectangular. \tcb{This requires a specific and new formulation of the charging problem, involving discrete action sets for the EVs in particular.} The considered cost function comprises of three components: the distribution transformer aging, the distribution \tcb{energy} losses, and a component inherent to the EV itself (e.g., the battery charging monetary cost). Charging start times are determined by the proposed distributed algorithm, whose analysis is conducted by using game-theoretic tools such as ordinal potential games. Convergence of the proposed algorithm is shown to be guaranteed for some important special cases. Remarkably, the performance loss w.r.t. the centralized solution is shown to be small. Simulations, based on realistic public data, allow one to gain further insights on the issues of convergence and optimality loss and provide clear messages about the tradeoff associated with the presence of the three components in the considered cost function. While simulations show that the proposed charging policy performs quite similarly to existing (continuous) charging policies such as valley-filling-type solutions when the non-EV demand forecast is perfect, they reveal an additional asset of rectangular profiles in presence of forecasting errors. 

\end{abstract}


\IEEEpeerreviewmaketitle

\section{Introduction}
\label{sec:introduction}


The deployment of electric vehicles (EVs) at a large scale is envisioned to have a significant impact on the existing and future energy networks \cite{ClementNyns10}. In the present paper, the impact on the grid is assessed in terms of residential distribution network (DN) costs. To be more specific, the main goal pursued is to optimize EV charging schedules to minimize a cost resulting from a linear combination of the residential distribution transformer aging and the distribution \tcb{energy} losses. \tcb{While the emphasis will be put on these two particular cost functions in this paper, all analytical results presented here apply to a large variety of costs; these results include the proposed problem formulation, the distributed charging algorithm, and its analysis. Among possible problems which can be considered we may mention voltage regulation \cite{beaude13}, harmonic distortion \cite{gomez03}, and reactive power management \cite{mojdehi14}. Concerning the component of the cost function which is referred as to the transformer component, the following has to be noted. Mathematically, it may correspond to any function of the past load levels; transformer aging is one possible instance and is the one made for the conducted numerical analysis. In practice,} as explained in \cite{Gong2011,McBee2009,Argade2012,Grahn2011,Rutherford2011,Hilshey2013,Turker12}, optimizing a long-term criterion such as the residential transformer lifetime becomes an important concern in the presence of EVs. Indeed, transformers might have to operate in a regime where aging is accelerated. In the European Union, about 5 millions of distribution transformers are used and about $70\%$ of transformer failures are due to aging \cite{Zhang2007}. Concerning \tcb{energy} losses in the DN, they represent the most important fraction of power losses in the whole electricity network; according to \cite{Dyevre2010}, in France two thirds of \tcb{energy} losses are due to the DN. Despite the importance of the problem\footnote{The problem directly concerns DN operators and car makers but may also concern the EV users since they might be charged in an indirect manner the extra costs induced by the impact of EV charging operations on the grid.}, the impact of EV charging on distribution transformer aging and \tcb{energy} losses has only been addressed in a relatively small number of papers. Among the relevant related works we may cite \cite{Gong2011}\cite{Argade2012}\cite{Rutherford2011}\cite{Shao11,6243548,6837513}. The dominant approach adopted, which is well illustrated by \cite{Gong2011}, consists in exploiting a suitable model for the aging or \tcb{energy} losses, and assessing the impact of charging for simple \textit{scenarios}; for instance, two possible scenarios are that all EVs start charging at a given time of the day (e.g., at $7$ pm) or at random times. The \textit{algorithmic} aspect is however not developed. This is precisely what the present work proposes. 

The algorithmic aspect of the charging problem has been tackled in the literature but, mostly, for minimizing the monetary cost the user has to pay for recharging his vehicle. In this respect \cite{Gan13}\cite{Ma13} constitute relevant works. Additionally, in the present paper the focus is on  \textit{distributed} charging in the sense that we assume the existence of several decision-makers (DMs) and each of them has only partial control of the variables which affect the cost or payoff function under consideration. Distributed charging is relevant in at least two key scenarios: the scenario in which charging policies have to be computed by a single entity (e.g., an aggregator \cite{Wu2012} or a transformer computing device) but for complexity issues it may be required to optimize the variables separately, forming a set of virtual DMs; the scenario where each EV controls its own charging policy, meaning that there are effectively several physical DMs instead of a single one. In the latter scenario, as far as advanced charging policies are concerned, DMs will be automata embarked on the EVs. For a naive policy such as just deciding when to plug the EV to the grid, the DMs might be the EV users themselves but this is not the standpoint adopted in this paper. 

It turns out that \textit{game theory} is very well suited\footnote{In particular, as explained further, the convergence and efficiency analyses for the proposed algorithm are conducted by introducing a charging game in which each player has his own payoff function.} to address the distributed charging problem whether players or DMs are persons or machines. In this respect, game-theoretic tools have been applied to smart grids quite recently (see e.g., \cite{Saad2012} for a survey); very useful contributions include \cite{Wu2011,Wu2012,Agarwal2011,Mohsenian-Rad2010,Vytelingum2010}. In \cite{Wu2011}, these tools are used for the DN frequency regulation problem in the context of the interaction between EVs and an aggregator. Therein, the authors show the usefulness of a well-chosen pricing policy to incite users to charge their vehicle in order to regulate the frequency of the electrical system. In \cite{Wu2012}, game theory is exploited to design a coordination mechanism for the wind power integration. References \cite{Agarwal2011} and \cite{Mohsenian-Rad2010} use a similar methodology to study the more general problem of load balancing whereas \cite{Vytelingum2010} applies this methodology for micro-storage management in smart grids. 

Compared to the application-oriented works where game theory is used to optimize energy consumption at the user side (at home, by the EVs, etc.), the present work possesses several distinguishing features. Two of them are as follows. First, to the best of our knowledge, our work is the first to propose a distributed charging algorithm which can not only minimize an individual cost, which is inherent to the EV and only depends on its actions (the battery charging monetary cost typically), but also the DN costs; the latter are given by the residential distribution transformer aging and \tcb{energy} losses over the DN. Second, we want to know to what extent using rectangular charging profiles is relevant both in terms of implementation and performance. The present work therefore adopts a complementary approach to existing works on charging algorithms which typically assume continuous charging power levels (see e.g.,  \cite{Gan13}\cite{Ma13}\cite{Shinwari2012}); the solutions used in the latter references are considered in the simulations section for comparison purposes. Although the focus of this paper is on the EV charging problem and DN costs, the obtained analytical and numerical results can be re-exploited for other problems in smart grids such as the problem of allocating or scheduling stored energy. 

The paper is structured as follows. Sec. \ref{sec:system-model} provides the model of the system and the network cost function components under consideration. Sec. \ref{sec:formulation-of-the-problem} describes how the distributed EV charging problem is formulated. Sec. \ref{sec:algo} provides the description of the proposed distributed charging algorithm and its main properties. Sec. \ref{sec:numerical-analysis} allows one to assess numerically the performance of the developed algorithm. The paper is concluded by Sec. \ref{sec:conclusion} .


\section{System modeling}
\label{sec:system-model}

\begin{figure}[tbh]\vspace{-3mm}
\begin{center}
\resizebox{9cm}{5cm}{
\psset{xunit=1cm,yunit=1cm}
\psset{arrowscale=2}
\begin{pspicture}(0,0)(26,15)
\psline[ArrowInside=->,linecolor=gray,linewidth=2pt](9,13.5)(9,12)
\pscircle[linecolor=gray,linewidth=2pt](9,9.5){1}
\pscircle[linecolor=gray,linewidth=2pt](9,11){1}
\psline[ArrowInside=->,linecolor=gray,linewidth=2pt](9,8.5)(9,7.5)
\psline[ArrowInside=->,linecolor=blue,linewidth=2pt](9,7.5)(6.91,5.21)
\psline[ArrowInside=->,linecolor=red,linewidth=2pt](9,7.5)(3,5.5)

\psline[ArrowInside=->,linecolor=blue,linewidth=2pt](9,7.5)(10.59,2.71)
\psline[ArrowInside=->,linecolor=red,linewidth=2pt](9,7.5)(15,3)

\psline[linecolor=red,linewidth=1.5pt](1,4.5)(3,4.5)
\psline[linecolor=red,linewidth=1.5pt](3,4.5)(3,6)
\psline[linecolor=red,linewidth=1.5pt](3,6)(1,6)
\psline[linecolor=red,linewidth=1.5pt](1,6)(1,4.5)
\psline[linecolor=red,linewidth=1.5pt](1,6)(2,7)
\psline[linecolor=red,linewidth=1.5pt](2,7)(3,6)
\psarc[linecolor=blue,linestyle=solid,linewidth=1.5pt](5.5,3.8){2}{30}{150}
\pscircle[linecolor=blue,linewidth=3pt](4.5,4.8){0.3}
\pscircle[linecolor=blue,linewidth=3pt](6.5,4.8){0.3}
\psline[linecolor=blue,linewidth=1.5pt](3.77,4.8)(4.25,4.8)
\psline[linecolor=blue,linewidth=1.5pt](6.8,4.8)(7.23,4.8)
\psline[linecolor=blue,linewidth=1.5pt](4.75,4.8)(6.25,4.8)

\psline[linecolor=red,linewidth=1.5pt](15,2)(17,2)
\psline[linecolor=red,linewidth=1.5pt](17,2)(17,3.5)
\psline[linecolor=red,linewidth=1.5pt](17,3.5)(15,3.5)
\psline[linecolor=red,linewidth=1.5pt](15,3.5)(15,2)
\psline[linecolor=red,linewidth=1.5pt](15,3.5)(16,4.5)
\psline[linecolor=red,linewidth=1.5pt](16,4.5)(17,3.5)
\psarc[linecolor=blue,linestyle=solid,linewidth=1.5pt](12,1.3){2}{30}{150}
\pscircle[linecolor=blue,linewidth=3pt](11,2.3){0.3}
\pscircle[linecolor=blue,linewidth=3pt](13,2.3){0.3}
\psline[linecolor=blue,linewidth=1.5pt](10.27,2.3)(10.7,2.3)
\psline[linecolor=blue,linewidth=1.5pt](13.3,2.3)(13.73,2.3)
\psline[linecolor=blue,linewidth=1.5pt](11.3,2.3)(12.7,2.3)
\pscircle*[linecolor=gray](9,7.5){0.15} 
\psline[linewidth=1pt]{->}(13,8)(25.5,8)
\psline[linewidth=1pt]{->}(13,8)(13,13.5)
\psline[linecolor=red,linewidth=1.5pt](13,10)(15,9)
\rput(13,7.5){{ \bf\huge 12 am}}
\rput(11.5,13){{\bf\huge Power}}
\psline[linecolor=black,linestyle=dashed,linewidth=.3pt](15,9)(15,8)
\psline[linecolor=black,linestyle=dashed,linewidth=.3pt](13,9)(15,9)
\rput(15,7.5){{\bf\huge 4 am}}
\psline[linecolor=red,linewidth=1.5pt](15,9)(17,12)
\psline[linecolor=black,linestyle=dashed,linewidth=.3pt](17,12)(17,8)
\rput(17,7.5){{\bf\huge 8 am}}
\psline[linecolor=red,linewidth=1.5pt](17,12)(19,12)
\psline[linecolor=black,linestyle=dashed,linewidth=.3pt](19,12)(19,8)
\rput(19,7.5){{\bf\huge 12 pm}}
\psline[linecolor=red,linewidth=1.5pt](19,12)(21,11)
\psline[linecolor=black,linestyle=dashed,linewidth=.3pt](21,11)(21,8)
\rput(21,7.5){{\bf\huge 4 pm}}
\psline[linecolor=red,linewidth=1.5pt](21,11)(22.5,13)
\psline[linecolor=black,linestyle=dashed,linewidth=.3pt](22.5,13)(22.5,8)
\psline[linecolor=black,linestyle=dashed,linewidth=.3pt](13,13)(22.5,13)
\rput(22.8,7.5){{\bf\huge 7 pm}}
\psline[linecolor=red,linewidth=1.5pt](22.5,13)(23,13)
  \psline[linecolor=red,linewidth=1.5pt](23,13)(24,11.1)
 \psline[linecolor=red,linewidth=1.5pt](24,11.1)(24.5,11.8)
 \psline[linecolor=red,linewidth=1.5pt](24.5,11.8)(25,10.5)
 \psline[linecolor=black,linestyle=dashed,linewidth=.3pt](25,10.5)(25,8)
 \rput(25,7.5){{\bf\huge 12 am}}
\end{pspicture}
}\vspace{-6mm}\caption{Considered network topology. A typical exogenous (non-EV) demand profile $L^{\mathrm{exo}}_1,..., L^{\mathrm{exo}}_T $ is represented.}\label{fig:DN} \vspace{-7mm}
\end{center}
\end{figure}
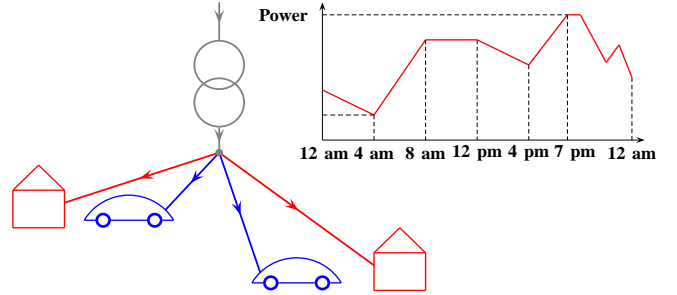

The goal of this section is to express the transformer aging and \tcb{energy} losses as a function of the sequence of the total demand power or load levels. To this end, the considered topology for the DN is as follows (see Fig.~\ref{fig:DN}). It consists of one residential transformer to which two groups of devices are connected: a set of EVs and a set of other electrical devices. The latter are assumed to induce a power demand which is independent of the charging policies and therefore called exogenous demand. The corresponding load is denoted by $L^{\mathrm{exo}}_t$, which is a deterministic function of the time which is assumed to be slotted and whose index is denoted by $t \in \mathcal{T}=\left\{1,...,T\right\}$; this function is always assumed to be known (except in the simulation part -Sec. \ref{sec:numerical-analysis}- where the influence of not forecasting it perfectly is assessed). The load induced by the EVs at time $t$ will be denoted by $L_t^{\mathrm{EV}}$ and the (average) total transformer load thus writes as $L_t = L^{\mathrm{exo}}_t + L_t^{\mathrm{EV}} $. Note that our analysis concerns the impact of EV charging in terms of active power; refinements related to the reactive power are left as a possible extension.


\subsection{Residential distribution transformer aging model}
\label{subsec:transfo-aging}

The most influential parameter for the transformer aging is known to be the hot-spot (HS) temperature \cite{IEEE1995}\cite{IEC}. Indeed, the transformer isolation damage is directly related to the HS temperature: aging is accelerated (decelerated) when the HS temperature is above (resp. below) its nominal value. In general, the HS temperature depends on the history of transformer load levels and that of the ambient temperature levels. Even though the HS temperature at time $t$ depends on the sequence of ambient temperature levels (see \cite{IEEE1995}\cite{IEC} for more details), we will not make this dependency explicit in our notations\footnote{The choice of a particular sequence of ambient temperature levels only intervenes in the numerical analysis.}. Indeed, this sequence cannot be controlled and in particular it does not depend on the power demand induced by the EVs. The HS temperature at time $t$ is given by

\begin{equation}
\label{ModelHSFonctInertia}
\theta^{\mathrm{HS}}_{t} = 
F^{\mathrm{HS}}_{t} \left( \boldsymbol{L}^t   \right) = F^{\mathrm{HS}}_{t} \left(L_1,...,L_t   \right)
\end{equation}
where $ \boldsymbol{L}^t=\left(L_1,...,L_t   \right)$ represents the sequence of total transformer load levels up to time $t$; the quantity $L_t$ is related to the charging power levels in Sec. \ref{sec:formulation-of-the-problem}. By default, no particular assumption will be made on the function $F^{\mathrm{HS}}_{t}$ throughout this paper. A few results will be provided for some special cases such as the memoryless case. Indeed, a special case of (\ref{ModelHSFonctInertia}) is when the HS temperature does not depend on the whole history of load levels but only on the current load level:
\begin{equation}
\label{ModelHSFonct}
\theta^{\mathrm{HS}}_{t} = \widetilde{F}^{\mathrm{HS}}_{t}\left({L}_{t} \right)  \textrm{.}
\end{equation}
Case (\ref{ModelHSFonctInertia}) (resp. (\ref{ModelHSFonct})) will be referred to as the case with (resp. without) thermal inertia. In \cite{IEEE1995,IEC,Rivera2008a} realistic models for the HS temperature evolution are provided, which allow one to have specific examples for  $F^{\mathrm{HS}}$ and $\widetilde{F}^{\mathrm{HS}}$. The corresponding specific expressions are only exploited in the numerical analysis. The transformer aging acceleration factor at time $t$, or aging for short, is assumed to be a function of the HS temperature at time $t$ and is denoted by $A_t$. Again, unless explicitly mentioned, no particular assumption will be made on the function $A_t$. A typical choice, which will be assumed in the numerical part, is as follows:
\begin{equation}
\label{EqAging}
A_t = e^{a \theta^{\mathrm{HS}}_{t} + b } = e^{a F^{\mathrm{HS}}_{t} \left( \boldsymbol{L}^t   \right)+ b } 
\end{equation}
where $a >0, b <0$ are some constants (see \cite{IEEE1995}\cite{IEC}). 

\tcb{\textit{Remark 1.}  All analytical results derived in this paper will only assume that $A_t$ is a function of the past load levels $\left(L_1(\mb{s}), ..., L_t(\mb{s})\right)$. On the other hand, numerical results will we based on the particular choice corresponding to (3).}

\subsection{Distribution network \tcb{energy} losses model}
\label{subsec:joule-losses}

Distribution \tcb{energy} losses mainly come from the transformer and the lines between the transformer and the active electrical devices. At the transformer level, both voltage and frequency are assumed to be fixed and thus independent of the load level induced by the EVs. This is why no-load losses will not be considered here. As the purpose of this paper is not to obtain a model for \tcb{energy} losses which is as advanced as possible, we assume a simple scenario which is sufficient to show to what extent the mathematical problem formulation and charging scheme are impacted. \tcb{Energy} losses on time-slot $t\in\mathcal{T}$ are assumed to express as
\begin{equation}
\label{JouleLosses}
J(L_{t}) = \left(R_{\textrm{transfo}}+R_{\textrm{line}}\right)  \left(  L^{\mathrm{exo}}_t + L_t^{\mathrm{EV}}  \right)  ^2
\end{equation}
where $R_{\textrm{transfo}}$ and $R_{\textrm{line}}$ are the transformer resistance and equivalent line resistance respectively; $L_{t} = L^{\mathrm{exo}}_t + L_t^{\mathrm{EV}} $ is the total transformer load on time-slot $t$. 
As mentioned previously, the assumed model can be improved but the retained model has at least three attractive features: 1) provided that continuous charging power levels are allowed, the minimization of $\sum_{t=1}^T J(L_t) $ w.r.t. $ (L_1^{\mathrm{EV}}, ...,  L_T^{\mathrm{EV}})$ corresponds to a valley-filling (VF) solution which is a well-known scheme \cite{Shinwari2012}. VF is performed over the sequence $(L^{\mathrm{exo}}_1,...,L^{\mathrm{exo}}_T)$. Here this sequence corresponds to the exogenous demand profile but it may also represent a sequence of prices (see e.g., \cite{Gan13});  2) the assumed \tcb{energy} losses model allows us to leave the load flow problem as a separate problem which might be handled with through an extension of our work \cite{Zhu12a}; 3) the model is relevant when line \tcb{energy} losses are dominated by transformer \tcb{energy} losses.
\section{Formulation of the EV charging problem }
\label{sec:formulation-of-the-problem}

This section aims at formulating in a distributed manner (as motivated in Sec. \ref{sec:introduction} the EV charging problem. The problem is said to be distributed because the variables which affect the payoff (or cost) functions of interest are not controlled jointly but separately. Since the EV charging profiles are imposed to be \textit{rectangular}, these variables correspond to the charging start times of the different EVs. The number of consecutive time instances or time-slots required to have the battery charged or to reach a required state of charge (SoC) for the next trip of EV $i$ is denoted by $C_i$ while the effective charging start time for EV $i$ is denoted by $s_i$. The individual payoff function which has to be maximized for EV $i \in \mathcal{I}$, $\mathcal{I} =\left\{1, ..., I\right\}$, is assumed to have the following form:
\begin{equation}
\label{PersonalCost}
u_{i}(s_{1},...,s_I) = -f_{i}\left(g^{\mathrm{DN}}_{i}\left(\bm{s}\right)
+g^{\mathrm{EV}}_{i}\left(s_{i}\right)\right)
\end{equation}
where $\boldsymbol{s} = (s_1, ..., s_I)$ and $f_i$ is the individual pricing function for
  user $i$, which is assumed to be strictly increasing; for instance, it may translate the technological costs induced by charging into a monetary cost. The function $g^{\mathrm{DN}}_{i}\left(\bm{s}\right)$ represents the cost associated with the DN. This cost is chosen to be a linear combination of the transformer aging and \tcb{energy} losses; its exact expression is provided a little further. The function $g^{\mathrm{EV}}_{i}$ can be any single-variable function of $s_i$. It is an individual cost which only concerns EV $i$. It may model the impact of the start time on EV $i$ battery aging, the individual electricity fare for user $i$, or its preference in terms of availability. For example, if the sequence of prices for EV $i$ over the time period of interest
is denoted by $(\pi_{i,1},..., \pi_{i,T})$ (i.e., the price is a function
of the time only), then a suitable choice for
$g_i^{\mathrm{EV}}$ might be
\begin{equation}
\label{eq:costMon}
g_i^{\mathrm{EV}}(s_i)= \beta \sum_{t=s_i}^{s_i+C_i-1} \pi_{i,t}
\end{equation}
where $\beta \geq 0$ is a weight which allows the tradeoff between individual EV preferences and DN costs to be tuned. The price model can even be more complicated mathematically than what is assumed in
\eqref{eq:costMon} e.g., by assuming that $\pi_{i,t}$ is a function of the current number of EVs charging at time $t$ i.e., $\pi_{i,t} = \phi_{i,t}(n_t)$ (see \cite{Ma13} for more details on this model). One of the main properties (namely, potentiality) of the charging game studied in Sec. \ref{sec:algo} is retained.
  
To explicit $g^{\mathrm{DN}}_{i}$ as a function of $\bm{s}$, some notations have to be introduced. Let $\widetilde{n}_{t}$ and $n_{t}$ respectively denote the numbers of EVs starting to charge and charging at time $t$. These quantities are related to $\boldsymbol{s} = (s_1, ..., s_I)$ by
\begin{equation}
\label{NEVBeginningCharge}
\widetilde{n}_{t}(\bm{s})=\displaystyle \sum_{i=1}^{I} \mathds{1}_{[s_{i}=t]} \ \text{and} \ 
n_{t}(\bm{s}) = \displaystyle \sum_{i=1}^{I} \displaystyle \sum_{t'=1}^{C_{i}} \mathds{1}_{{[s_{i}=t-C_{i}+t']}}
\end{equation}

where $\mathds{1}_{[.]}$ is the indicator function. To evaluate the impact of the charging policies on the DN, the key quantity is the total transformer load or consumed power. The $T-$dimensional sequence of total load levels is now also denoted by $\bm{L}^T(\bm{s})= \left(L_1(\bm{s}), ..., L_T(\bm{s}) \right)$ and the total load on time-slot $t$ is given by
\begin{equation}
\label{TotLoad}
L_t(\bm{s}) =  L^{\mathrm{exo}}_t + L_t^{\mathrm{EV}}(\bm{s}) = L_t^{\text{exo}} + P n_t(\bm{s})
\end{equation}
where $P$ is the common charging power of all the EVs. Using the introduced notations, it is now possible to express $ g^{\mathrm{DN}}_{i}$:
\begin{equation}
    \label{GridCost}
      g^{\mathrm{DN}}_{i}\left(\bm{s}\right) =
      \displaystyle \sum_{t \in \mathcal{W}_i(s_i)} \alpha A_{t}\left(\bm{L}^{t}\left(\bm{s}\right)\right)
      +(1-\alpha)J\left(L_{t}\left(\bm{s}\right)\right)
\end{equation}
where $0 \leq \alpha \leq 1$ is the weight given to transformer aging relatively to \tcb{energy} losses and $\mathcal{W}_i(s_i)$ is a discrete set which represents the time window over which the EV $i$ is considered to be influential on the DN cost. We will dedicate more attention to two special cases of practical interest for $\mathcal{W}_i(s_i)$. The first case is when $\mathcal{W}_i(s_i) =  \lbrace s_i,\cdots,s_i+C_i-1 \rbrace \triangleq \mc{W}_i$ which means that each EV individual payoff is only related to the period of time over which the corresponding EV is active. The second case corresponds to $\mathcal{W}_i(s_i)  = \mc{T}$. The parameter $\alpha$ can be seen as a simple way of tuning the tradeoff between a short-term cost (\tcb{energy} losses) and a long-term cost (transformer aging). 

Now let us specify the action set for each EV. The arrival and departure time of EV $i \in \mathcal{I}$ are denoted by $a_i \in \mathcal{T}$ and $d_i \in \mathcal{T}$ respectively. As it is assumed that an EV has to charge its battery within the total time window, the action set for EV $i$ is chosen to be
\begin{equation}
\label{DebCharge}
\mathcal{S}_{i}=\left\{a_{i},a_{i}+1,...,d_{i}-C_{i}+1\right\}.
\end{equation}
The EV action profile $\bm{s}$ therefore lies in $ \mathcal{S}=\prod_{i=1}^{I}  \mathcal{S}_{i}$; the standard notation $\bm{s}_{-i}=\left(s_{1},...,s_{i-1},s_{i+1},..,s_{I}\right)$, $I\geq 2$, will be used for referring to the reduced action profile in which user $i$'s action is removed. A special case of interest is when all users have the same charging constraint, that is $\forall i, \, a_{i}=a, d_{i}=d, C_{i}=C$. This case will be said to be symmetric. In the symmetric case, it will be assumed, without loss of generality, that $a=1$ and $d=T$. 

\tcb{\textit{Remark 2.}  Maximizing the sum of weighted individual costs provides a Pareto optimal point when the cost region is convex \cite{geoffrion68}. By changing the weights, one moves along the Pareto frontier of the cost region, which represents the best which can be done in terms of tradeoff between aging and energy losses. This therefore gives a motivation for considering the linear combination. A natural question is about whether the cost region is convex. In fact, the feasible \textbf{average} cost region is convex whenever the cost is averaged over a large number of stages (which is the number of days here). This follows by a time-sharing argument: if there exist two charging strategies which achieve each a given pair of individual costs, then any convex combination of these pairs can be achieved by using the two strategies with the appropriate fraction of the time.} 

\textit{Remark 3.} As already mentioned, charging profiles are assumed to be rectangular, which is why the charging power can only have two possible levels namely either $0$ or $P$. Among the motivations for considering rectangular charging profiles we may mention the following: 1) An important argument is that rectangular profiles are currently being used for existing EVs (e.g., the EVs built by the French car maker RENAULT) and not only in papers; 2) for a given charging start time, charging at full power without interruption minimizes the delay to charge; 3) one technological reason is that they allow one to manage the EV battery aging. 
Battery aging seems to be  accelerated when the charging operation comprises interruptions \cite{Gan2012}; 4) as shown further in Sec. \ref{sec:numerical-analysis}, rectangular profiles are also fully relevant for cost functions with memory such as the transformer aging. The charging start time turns out to be a very influential variable, confirming the observations made in \cite{Rutherford2011}; 5) more specifically,  from an optimal control theory perspective, rectangular profiles may be optimal. This
 happens for instance when the state (i.e., the hot-spot temperature)
 is monotonically increasing with the
 control (i.e., the charging power). If, for the time window of interest,
 the transformer temperature can only increase, it is optimal to
 delay the transformer heating. The optimal solution is then to
 start charging as late as possible i.e., to charge at maximal power
 at the end of the considered time window and charge at zero power
 before. The duration of the corresponding charging profile, which is rectangular, is given by the desired final state of charge; 6) the charging start time will be seen to be less sensitive to forecasting errors on the exogenous demand than VF-type solutions; 7) profiles without interruption are even required in some important scenarios encountered with home energy management \cite{Bapat11}\cite{Chavali14}.

\section{A new distributed charging algorithm}
\label{sec:algo}

\subsection{Motivations}

As explained in Sec. \ref{sec:introduction}, considering distributed charging policies is relevant in at least two key scenarios. Assume a scenario (say \textit{Scenario} $c$) in which the charging policies are computed by a single decision-making entity (e.g., an aggregator or a transformer computing device) and the maximization of a quantity such as the sum-payoff function  $\sum_{i=1}^I u_i(\boldsymbol{s})$ is pursued. Finding an optimal solution may largely exceed the available computational capacity; an exhaustive search would roughly involve $T^I$ tests. Therefore, even if there is one single decision-making entity, it may be required to optimize the variables of $\boldsymbol{s}$ separately. Now, if we assume a scenario (say \textit{Scenario} $d$) in which each EV controls its own charging policy, there are $I$ physical DMs and the problem is distributed by nature. In both scenarios, the total computational complexity of the proposed algorithm will be seen to be typically linear in the product $T \times I$, showing the dramatic reduction in terms of complexity allowed by the used distributed implementation.

The proposed distributed algorithm is based on a procedure which is called the sequential best-response dynamics (BRD) in game theory literature (see e.g., \cite{Fudenberg1991}\cite{lasaulce-book-2011}). One of the strong motivations for selecting such a procedure for the problem under investigation is that convergence of the associated iterative algorithm can be guaranteed with overwhelming probability. Additionally, convergence is very fast, which is useful both in Scenario $c$ to avoid unnecessary computations and in Scenario $d$ to manage the amount of signaling between the EVs and the aggregator. Other arguments in favor of using the BRD will be provided  further.

To clearly indicate that no strategic assumption such as rationality or complete information is required on the DMs which implement the BRD, the description of the algorithm has been separated from its analysis. The analysis relies on the use of game-theoretic tools such as the powerful notion of potentiality, which guarantees the existence of a Nash point in a game and the convergence of the BRD to a Nash point. It turns out that the charging game under consideration is effectively potential and therefore makes the BRD a good candidate for computing the charging start instants.

\subsection{Description of the algorithm}

The proposed algorithm to determine the vector of start times $(s_1,...,s_I)$ is an iterative algorithm which is inspired from the sequential BRD. The algorithm is performed offline, which means that the decisions which intervene in the algorithm are intentions but not decisions which have effectively been taken; only the decisions obtained after convergence will be effective and implemented online. Once the charging instants are computed, the EV can effectively charge their battery according to the schedule determined. In its most used form, the BRD operates sequentially such that DMs update their strategies in a round-robin manner. Within round $m$ (with $m\geq1$) the action chosen by DM $i$ (which can be virtual or physical depending on the assumed scenario) is computed as (\ref{eq:update}).
The proposed procedure is translated in pseudo-code through Algorithm 1.\\
\begin{algorithm}
\emph{Initialize the round index as $m=0$. Initialize the vector of charging start times as $\bm{s}^{(0)}$.} \\
\While{$\left\| \bm{s}^{(m)} - \bm{s}^{(m-1)}\right\|  >  \delta \hspace{0.1cm}$ and $\hspace{0.1cm} m \leq M \hspace{0.1cm}$}{
\emph{\textbf{Outer loop}. Iterate on the round robin phase index: $m=m+1$. Set $i=0$. \\}
\emph{\textbf{Inner loop}. Iterate on the DM index: $i=i+1$. Do: \\}
\begin{align}\label{eq:update}
s_i^{(m)} \in \arg \max_{s_i \in \mc{S}_i} u_{i}(s_1^{(m)}, s_2^{(m)},..., s_i, \nonumber \\
s_{i+1}^{(m-1)}, ..., s_I^{(m-1)})
\end{align}
where $s_i^{(m)}$ stands for action of DM $i$ in the round robin phase $m$. Stop when $i=I$ and go to \textbf{Outer loop}.}
\caption{\label{algo} The proposed distributed EV charging algorithm.}
\end{algorithm}
\\
\textit{Comments on Algorithm 1.}\\
$\bullet$ In (\ref{eq:update}), when the argmax set is not a singleton, $s_i^{(m)}$ is randomly drawn among the maximum points.\\
$\bullet$ The quantity $\delta \geq 0$ in Algorithm 1 corresponds to the accuracy level wanted for the stopping criteria in terms of convergence.\\
$\bullet$ To update the charging power levels $m$ times, $m \times I$ iterations are required.\\
 $\bullet$ The order in which DMs update their action does not matter to obtain convergence (see e.g., \cite{Bertsekas-book-95}). However, simulations which are not provided here indicate that some gain in terms of convergence time can be obtained by choosing the order properly. A good rule seems to be to start updating at iteration $m$ the EV decisions in an increasing order in terms of start times as obtained per iteration $m-1$, which makes the order iteration-dependent.\\
$\bullet$ The knowledge required to implement Algorithm 1 is scenario-dependent. In Scenario $c$ in which each decision is computed by a single entity (the transformer typically), the vector of effective charging instants can be computed from its initial value $\boldsymbol{s}^{(0)}$, the (forecasted) sequence of exogenous loads $(L_1^{\mathrm{exo}},...,  L_T^{\mathrm{exo}})$, and the parameters which intervene in the payoff functions; the latter include in particular the EV mobility data $(a_i)_{i \in \mc{I}}$, $(d_i)_{i \in \mc{I}}$, and $(C_i)_{i \in \mc{I}}$. In Scenario $d$ in which the EVs themselves update their decision, messages have necessarily to be exchanged between the transformer and the EVs. A possible communication protocol is as follows. Without knowing anything about the exogenous demand $\bs{L}^{\mathrm{exo}}$ and the one associated with the other EVs, EV automaton $i$ chooses a start time say $s_i^{(0)}$ and reports this to the transformer. The latter aggregates the received signals and replies to the EVs by sending them the predicted sequence of total load levels $L_t^{(0)} =  L_t^{\mathrm{exo}} + P n_t(\boldsymbol{s}^{(0)})$, $t\in \mc{T}$. Therefore, EV $1$ updates its intended start time as $s_1^{(1)} \in \arg \max_{s_1} u_1(s_1, \bs{s}_{-1}^{(0)} )$ and reports this change to the transformer: the latter updates the aggregate signal into $L_t^{(1)} =  L_t^{\mathrm{exo}} + P n_t(\bs{s}^{(1)})$, $t\in \mc{T}$, where $\bs{s}^{(1)} = (s_1^{(1)}, s_2^{(0)},..., s_I^{(0)})$. This signal is sent to all the EVs but only EV $2$ is able to update its intended start time. When all EVs have updated their start time at least once, a new updating round can start. In practice, it might happen that much less knowledge is available. For instance, if only the knowledge of the forecast exogenous demand $L^{\mathrm{exo}}$ is available then it is always possible to apply Algorithm 1 on the cost function (\ref{PersonalCost}) with $\alpha =0$. If the corresponding charging scheme is used, it will induce, in general, a certain loss of optimality.\\ 
$\bullet$ A variation of Algorithm 1 can be obtained by updating the charging policies simultaneously. The main reason why we have not considered the parallel version is that it is known that there is no general analytical result for guaranteeing convergence \cite{lasaulce-book-2011}. When converging, the parallel implementation is faster but since start times are computed offline, convergence time may be seen as a secondary feature. 

\subsection{Convergence analysis}
\label{sec:sub-convergence-analysis}

One of the powerful links between distributed optimization and game theory is that scenarios involving several individual optimizers or DMs which update their decisions over time may converge to a Nash equilibrium (NE) of a certain game. This is one of the reasons why we now define the game of interest, that we will refer to as the charging game. The main purpose of this section is to show that, under additional realistic assumptions, this game is an ordinal potential game (OPG) \cite{Monderer1996}. The intuition behind this is  that the considered cost functions all depend on the EV charging loads through the sum-load. This type of structures, which is present in games which are called aggregate games (a congestion game \cite{rosenthal73} is a special instance of them), may lead to the existence of a potential function. Although exact potentiality is typically not available here, it turns out that ordinal potentiality is available under some mild conditions. In particular, the latter property guarantees the convergence of the distributed charging algorithm proposed in Sec. \ref{sec:algo}. A game under strategic-form is given by an ordered triplet which respectively comprises the set of DMs or players, their strategy sets, and their payoff functions (see e.g., \cite{lasaulce-book-2011}). The charging game of interest is defined as follows.

\begin{definition}[Charging game] The charging game is the triplet $\mc{G} = \left( \mc{I}, \left(\mc{S}_i\right)_{i\in\mc{I}}, \left(u_i\right)_{i\in\mc{I}} \right)$ whose elements are defined by (\ref{PersonalCost}) and (\ref{DebCharge}).
%
%
\end{definition}
Interestingly, this game is an OPG for typical scenarios encountered in practice. An OPG is defined as follows \cite{Monderer1996}.
\begin{definition}\label{def:OPG}[Ordinal Potential Game] A game whose payoff functions are $\left(u_i\right)_{i\in\mc{I}}$ is an OPG if there exists a function $\Phi$ such that $\forall i \in \mathcal{I}, \, \forall \bm{s}=(s_{i},\bm{s}_{-i}), \, \forall s_{i}^{\prime} \in \mc{S}_i$,
\begin{equation}
\label{GameOrdPotential}
u_{i}(s_{i}^{\prime},\bm{s_{-i}}) \geq u_{i}(s_{i},\bm{s_{-i}}) \Leftrightarrow  \Phi(s_{i}^{\prime},\bm{s_{-i}}) \geq \Phi(s_{i},\bm{s_{-i}}) \textrm{ .} \nonumber
\end{equation}
\end{definition}
Since the function $\Phi$ does not depend on the player index, the game analysis amounts, to a large extent, to analyzing an optimization problem. This attractive property is available at least in the scenarios defined by the next proposition.

\begin{proposition}[Potential property of the charging game]\label{PropPotGame} If one of the following conditions is met, then $\mc{G}$ is an OPG: (a) $\forall i \in \mc{I}, \ \mathcal{W}_i(s_i) = \mathcal{W}$, where $\mathcal{W}$ is any discrete set which is independent of the player index $i$ and $s_i$; (b) \tcb{$\forall i \in \mc{I}, \ \mathcal{W}_i(s_i) = \lbrace s_i,\cdots,s_i+C_i-1\rbrace$ and} there is no thermal inertia i.e., $\theta^{\mathrm{HS}}_{t}$ is given by (\ref{ModelHSFonct}). Two potential function candidates are respectively given by: (a) $\Phi_{(a)}\left(\bm{s}\right)   =  -\displaystyle \sum_{t \in \mathcal{W}} \alpha A_{t}\left(\bm{L}^t\left(\bm{s}\right)\right)
      +(1-\alpha)J\left(L_{t}\left(\bm{s}\right)\right)  -\displaystyle \sum_{i\in \mathcal{I}}g^{\mathrm{EV}}_{i}\left(s_{i}\right)$; (b) $\Phi_{(b)}\left(\bm{s}\right)   =  -\displaystyle \sum_{t \in \mathcal{T}} \displaystyle \sum_{v_{t}=0}^{n_{t}(\bm{s})}  \alpha  A_{t}\left(L_{t}(v_t) \right) + (1-\alpha)J\left(L_{t}(v_t) \right)  -\displaystyle \sum_{i\in \mathcal{I}}g^{\mathrm{EV}}_{i}\left(s_{i}\right)$ where $L_t(v_t) = L_t^{\mathrm{exo}} + P v_t$.
\end{proposition} 
The proof is provided in App. \ref{app:proofPot}. The two scenarios in which the game is potential are clearly of practical interest. Assuming all the EVs to have a common time window $\mc{W}_i(s_i) =  \mc{W}$ for charging (Scenario (a)) means that the total cost associated with the DN is accounted for by all the EVs. This is clearly the most interesting scenario for the aggregator. However, from the user's standpoint, this may be considered as unfair or not acceptable. In such a case, it is more realistic that user $i$ be only charged a cost which corresponds to the period over which his battery is effectively recharged namely, in the time interval defined by $\mc{W}_i(s_i) = \{s_i,s_i +1, ..., s_i+C_i-1 \}$ (note that an alternative way of individualizing the cost would be to use a cost sharing policy such as in \cite{Wu2012} which would lead us to use a common window $\mc{W}$ but weighting the total cost by $\frac{C_i}{\sum_j C_j}$ for user $i$). This general scenario is mathematically more involving than Scenario (a). It turns out that it becomes quite simple to be analyzed for transformers with low thermal inertia, which leads to Scenario (b). Obviously, when \tcb{energy} losses represent the dominant cost ($\alpha \rightarrow 0$), the game is always potential. All these comments lead us to the next proposition which is the key result of this section.

\begin{proposition}[Sufficient conditions for convergence]\label{CVAlgo} Algorithm~\ref{algo} converges if one of the following conditions is met: (a); (b); $\alpha = 0$.
\end{proposition}

As proved in \cite{Monderer1996}, the sequential BRD converges in ordinal potential games. The proof therefore follows. When the charging game does not meet none of the conditions above, a deeper analysis has to be conducted. In this paper, our choice is to exploit Monte-Carlo simulations to provide additional insights on the convergence issue. In Sec. \ref{sec:numerical-analysis}, the empirical convergence probability is assessed for other typical scenarios and is shown to be high. When converging, the proposed algorithm converges to an NE of the charging game $\mc{G}$. The motivation of the next section is to analyze the existence, uniqueness, and global efficiency of the convergence point(s) which are NE.

\subsection{Equilibrium analysis}
\label{sec:equilibrium-analysis}

A pure NE  of $\mc{G}$ is a point which meets a certain condition of stability \cite{Nash-1950}. Formally, it is defined as follows.

\begin{definition}[Pure NE] The action profile $\bm{s}^{*}=\left(s_{1}^{*}, ...,  s_{I}^{*}  \right) \in \mathcal{S}$ is a pure NE if $\ \forall i, \forall s_{i}, \ \ u_{i}(s_{i},\bm{s_{-i}}^{*}) \leq u_{i}(\bm{s}^{*})$.
\end{definition}

The \textit{existence} of a pure NE in $\mc{G}$ is ensured under the conditions assumed in the following proposition.

\begin{proposition}[Existence of a pure NE]\label{PropExistPureNE} In Scenarios (a), (b), or $\alpha=0$, the game $\mc{G}$ has at least one pure NE.
\end{proposition}

The proof of this result follows from the fact that in any of the scenarios of Prop. \ref{PropPotGame}, $\mc{G}$ is an OPG and after \cite{Monderer1996} the existence of a pure NE is guaranteed. On the other hand, and this is common in discrete games, \textit{uniqueness} does not hold. To prove this set $T=5$, $I=3$, $C_{i}=C=2$, $a_{i}=a=1$, $d_{i}=d=5$, $\bm{L}^{\mathrm{exo}}=(1,2,3,2,1)$, and $P=1$. For a small transformer inertia i.e., a small $T^0$ (see Sec. \ref{sec:numerical-analysis}), it can be checked that  $\bm{s}^{*}=(1,1,4)$, $\bm{s}^{**}=(1,4,1)$, and $\bm{s}^{***}=(4,1,1)$ are NE, showing that uniqueness is not guaranteed in general. Since NE uniqueness is not guaranteed in general, measuring the \textit{efficiency} of the worst NE is important. The most usual way of assessing the impact of decentralization in a non-cooperative game has been formalized in \cite{Papadimitriou2001} by defining the notion of price of anarchy. Rather, we will slightly modify the latter notion as the price of decentralization ($\mathrm{PoD}$), which we define below; the merit of the proposed definition is just that the price is effectively zero when a distributed algorithm or procedure leads to an equilibrium point which performs as well as the centralized solution in terms of \textit{sum-payoff}  $w = \sum_{i \in \mathcal{I}} u_i$.


\begin{definition}[$\mathrm{PoD}$]
The $\mathrm{PoD}$ of $\mc{G}$ is defined by
\begin{equation}
\mathrm{PoD}= 1-  \frac{\ds{\max_{\mb{s} \in \mc{S}}} \ \ w(\mb{s})}{   
\ds{\min_{\mb{s} \in \mc{S}^{\mathrm{NE}}} w(\mb{s})}}
\end{equation}
where $\mc{S}^{\mathrm{NE}}$ is the set of NE of the game. 
\end{definition}
\tcb{It can be seen that $0 \leq \mr{PoD} \leq 1$ and the larger the $\mr{PoD}$, the larger the loss due to decentralization}. It is generally difficult to express the above quantity as a function of the game parameters \cite{lasaulce-book-2011}. This explains why these quantities are often more relevant from the numerical point of view. Nonetheless, it is possible to characterize it in some special cases. One of the cases where $\mathrm{PoD}$ can be characterized is the limit case of a large number of EVs, that is $I \rightarrow \infty$, having the same charging constraint $a_i=a$, $d_i=d$ and $C_i=C$, and the transformer has no thermal inertia. In this asymptotic regime, $\frac{n_t}{I}\rightarrow x_t \in \mathbb{R}$ represents the proportion of EVs charging at time $t$ and the analysis of the game $\mc{G}$ amounts to analyzing the so called \textit{non-atomic} counterpart $\mc{G}^{\textrm{NA}}$ of $\mc{G}$. In the latter game, the set of players is continuous and given by $\mathcal{I}^{\text{NA}} =\left[0,1\right]$. The action set of the EVs, $\mathcal{S}$, is defined as in (\ref{DebCharge}). In the regime of large numbers of EVs, the transformer load becomes $\bm{L}(\bm{x})=\bm{L^{\text{exo}}}+p\bm{x}$; the parameter $p$ is introduced in order for the exogenous demand to scale with $I$. Indeed, after (\ref{TotLoad}), when $I\rightarrow +\infty$, if kept fixed, the exogenous demand $L_t^{\mathrm{exo}}$ tends to vanish in comparison to the load induced by the EVs. This is the reason why we introduce the parameter $p$ (instead of $P$). The obtained non-atomic charging game can be proved to be an OPG and the following result concerning efficiency can be obtained.

\begin{proposition}[$\mathrm{PoD}$ in the non-atomic case ($I \rightarrow \infty$)]\label{PropPoA} Assume that: $a_i=a$, $d_i=d$, $C_i=C$, and $\mc{W}_i(s_i) = \mc{T}$; $g_i^{\mathrm{EV}}=0$; $L_t^{\mathrm{exo}}$ is a non-increasing (resp. non-decreasing) function of
$t$ on $\lbrace 1,\cdots,C\rbrace$ (resp. $\lbrace T-C+1,\cdots,T\rbrace$); Scenario (b) or $\alpha=0$ is considered. Then we have that $\mathrm{PoD} = 0$.
\end{proposition}

The proof of this result is provided in App. \ref{app:proofPoD}. This result has the merit to exhibit a scenario where decentralizing the charging decisions induces no cost in terms of global optimality for the sum-payoff. Note that, in particular, if the exogenous demand is either constant or negligible w.r.t. the demand associated with the set of EVs, the above assumption holds and there is therefore no efficiency loss due to decentralization. 

\textit{Remark 2.} The fact that there exist multiple Nash equilibria might be seen as a crucial point since one does not not know to which point Algorithm 1 will converge. However, Prop. IV.4 and the numerical results provided in the next section show that the PoD is typically small. This means that the worst Nash equilibrium and the best Nash equilibrium necessarily perform similarly, showing that equilibrium selection is not a crucial issue for the problem under consideration.

\section{Numerical analysis}
\label{sec:numerical-analysis}

We first provide the general simulation setup assumed by default while particular choices will be specified in the figure captions. Data corresponding to non-EV demand (or exogenous) profiles and the ambient temperature are taken from the ERDF French DN operator data basis. They concern France for the year $2012$ and can be found in \cite{erdf,Tamb}. Unless specified otherwise, simulations are performed over a year; the chosen time unit corresponds to $30$ min. We consider a $20$ kV/$410$ V transformer whose apparent power is $100$ kVA and nominal (active) power is $90$ kW. The transformer HS temperature evolution law is assumed to follow the ANSI/IEEE linearized Clause 7 top-oil-rise model, which is described in \cite{Rivera2008a}. The transformer lifetime is inversely proportional to the average aging: 
\vspace{-2mm}
\begin{equation}
\text{lifetime} = 40 \times T_{\mathrm{year}} \times \left(\ds{\sum_{t=1}^{T_{\mathrm{year}}}  A_t} \right)^{-1} \ (\text{years}) 
\end{equation}
where the non-EV or exogenous demand is normalized such that without EV $\text{lifetime} = 40$ years; here $T_{\mathrm{year}} = 366 \times T$. The instantaneous aging $A_t$ (see (\ref{EqAging})) is computed by choosing $a = 0.12$\degre C$^{-1}$, $b= -11$. The function $F^{\mathrm{HS}}_{t}$ is not described here but can be found in \cite{Rivera2008a}. To make the simulations reproducible we provide the values of the different parameters of $F^{\mathrm{HS}}_{t}$: $\Delta_t=0.5$ h; $T^{0}=2.5$h (thermal inertia) for all simulations concerning the transformer, unless specified otherwise (in which case we have that $T^{0}=0.5$h) ; $\gamma= 0.83$; $R=5.5$; $\Delta \theta^{\mathrm{O}}_{\textrm{FL}} = 55$\degre C; $\Delta \theta^{\mathrm{HS}}_{\textrm{FL}} = 23$\degre C; $q = 1$; $r=1$; $\theta^{HS}_{0}=98$\degre C. \tcb{Energy} losses are evaluated by choosing $R_{\textrm{transfo}} = R_{\textrm{line}} = 0.03 \ \Omega$. The load induced by one EV is $P=3$ kW. For each day, charging operations have to take place within the time window from $5$ pm (day number $j$) to $8$ am (day number $j+1$), which corresponds to $\mc{T}= \left\{1,2,...,30\right\}$. Concerning the EV mobility data, two scenarios will be considered. In \textit{Scenario (s)}, all EVs need $C_i=C=16$ time-slots of $30$ minutes each to completely recharge their $24-$kWh battery and $a_i=a=1$ ($5$ pm), $d_i=d=30$ ($8$ am); this scenario can be seen as the worst case. In \textit{Scenario (t)}, the mobility data are deduced from statistics taken from the French survey ENTD 2008 available in \cite{entd}: $a_i$, $d_i$, and $C_i$ are taken to be the closest integers of realizations of Gaussian random variables $\widetilde{a}_i \sim \mc{N}(4,1.5)$, $\widetilde{d}_i \sim \mc{N}(29,0.75)$ and $\widetilde{C}_i \sim \mc{N}(5.99,1.14)$. By default, Scenario (s) will be assumed. Finally, unless specified otherwise, the cost functions in \eqref{PersonalCost} are defined by: for all $i$, $g^{\textrm{EV}}_{i}=0$, which allows us to isolate the effects of the exogenous demand; $\mathcal{W}_i(s_i)=\left\{s_i,s_i+1,...,s_i+C_i-1\right\}$, and $f_i=\mathrm{Id}$. The plug-and-charge (PaC) policy is obtained by assuming that EVs start charging as soon as they plug to the grid according to the data of \cite{entd}. To assess the impact of not being able to forecast the non-EV demand perfectly we assume that the available non-EV demand profile is given by $\widetilde{L}_t^{\mathrm{exo}} = L_t^{\mathrm{exo}} + Z$, where $Z \sim \mathcal{N}(0,\sigma_{\mathrm{day}}^2)$. We define the forecasting signal-to-noise ratio by
\begin{equation}
 \mathrm{FSNR} = 10\log_{10}\left(\frac{1}{\sigma_{\mathrm{day}}^2} \times \frac{1}{T_{\mathrm{day}}} \sum_{t=1}^{T_{\mathrm{day}}} (L_t^{\mathrm{exo}})^2\right) (\mathrm{dB})
\end{equation} 
where $T_{\mathrm{day}} =48$. 

\begin{figure}
\begin{center}
\includegraphics[scale=0.62]{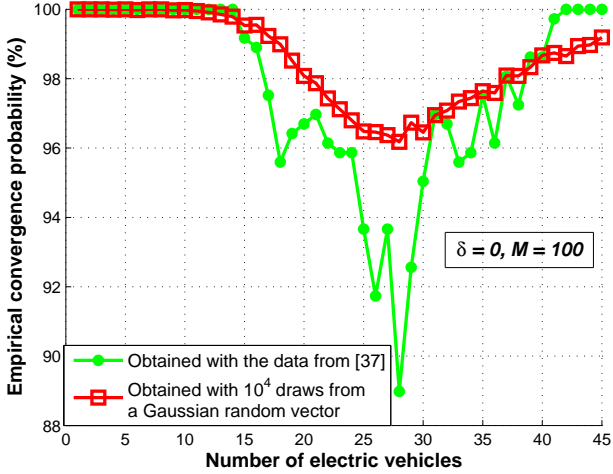}
\caption{Even if the assumption of Prop. \ref{CVAlgo} is not met, Algorithm 1 still converges with high probability. When the number of EVs ($I$) is low or large, the (empirical) probability of convergence tends to $1$. When the number of EVs takes intermediate values, the probability is typically above $90\%$.}
\label{fig:CVBRDyn}
\vspace{-7mm}
\end{center}
\end{figure}

\textit{Numerical convergence analysis.} In Sec. \ref{sec:sub-convergence-analysis}, we have provided sufficient conditions under which convergence is guaranteed. Here, we consider a scenario in which these conditions are not met (in fact the worst, with $\alpha=1$) and assess the probability of convergence of Algorithm 1. For Scenario (s), Fig. \ref{fig:CVBRDyn} represents the empirical probability of convergence against the number of EVs for the $366$ exogenous demand profiles from \cite{erdf} and for  $10 \ 000$ draws from a Gaussian random vector; the covariance matrix of the latter is taken to be $\sigma^2 \times \bm{I}_{\lbrace T \times T\rbrace}$, where $\sigma = 26$kW is estimated by using the data from \cite{erdf}. All simulations performed showed that only a few iterations are needed to obtain convergence, which is a quite typical behavior for sequential BRD-type iterative procedures \cite{lasaulce-book-2011}. 

\begin{figure}
\begin{center}
\includegraphics[scale=0.62]{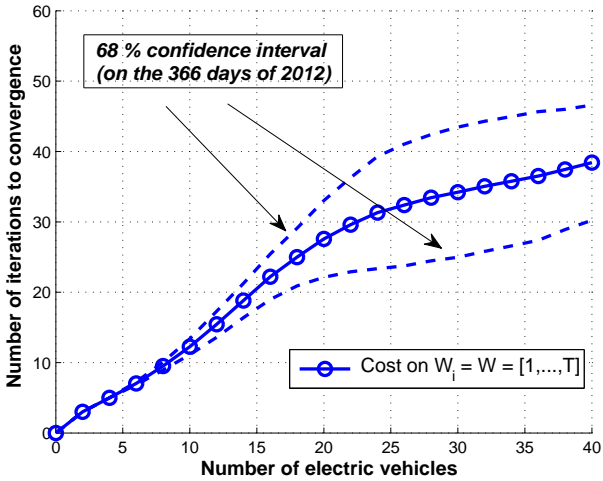}
\caption{Total number of iterations  of Algorithm 1  against the number of EVs ($I$). The dashed curves represent the $68\%$ confidence intervals. It is seen that Algorithm 1 is scalable in terms of convergence time. Here convergence time scales well with $I$ since variations are smooth.}  
\label{fig:vitCVbrd}
\vspace{-5mm}
\end{center}
\end{figure}

To elaborate further on this point, we provide Fig. \ref{fig:vitCVbrd}. It represents the total number of iterations needed for convergence of Algorithm 1 as a function of the number of EVs when $\mathcal{W}_i(s_i) = \mc{T}$. The middle curve corresponds to an average over the $366$ days of 2012 while the two others are associated with a chosen confidence interval of $68\%$. The total convergence time is observed to scale very well with the number of EVs indicating no scalability issues for convergence. Note that the scaling law seems almost piecewise linear here but many other simulations performed for more diverse scenarios show it might be "less linear" but always involves smooth variations.

\begin{figure}
\begin{center}
\includegraphics[scale=0.62]{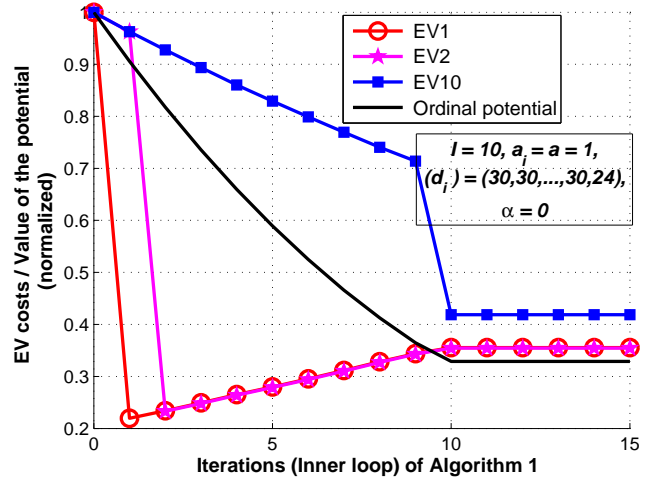}
\caption{For $10$ EVs, only $10$ iterations are needed to converge, which means that only one update per EV is sufficient to reach convergence, showing a very fast convergence. The figure also illustrates the fact that the potential is a Lyapunov function of the considered dynamical procedure.} 
\label{fig:CVbrd}
\vspace{-7mm}
\end{center}
\end{figure}

To conclude on the convergence issue, we consider the variations of the individual payoff functions of the EVs and that of a potential function over iterations; this is the purpose of Fig. \ref{fig:CVbrd}. For $I=10$, $\forall i, \ a_i=1$, $d_1=...=d_9=30, d_{10}=24$, the night of January 1, 2012 \cite{erdf}, and now $\alpha=0$ (only \tcb{energy} losses are taken into account), the figure shows that the potential always decreases over iterations, which illustrates the fact that it is a Lyapunov function of the considered dynamical procedure. As a very positive result, it is seen here that only one update per EV is needed to reach convergence. This very fast convergence behavior is typical when a memoryless cost is considered \cite{lasaulce-book-2011}, which is the case here with $\alpha=0$.

\begin{figure}
\begin{center}
\includegraphics[scale=0.65]{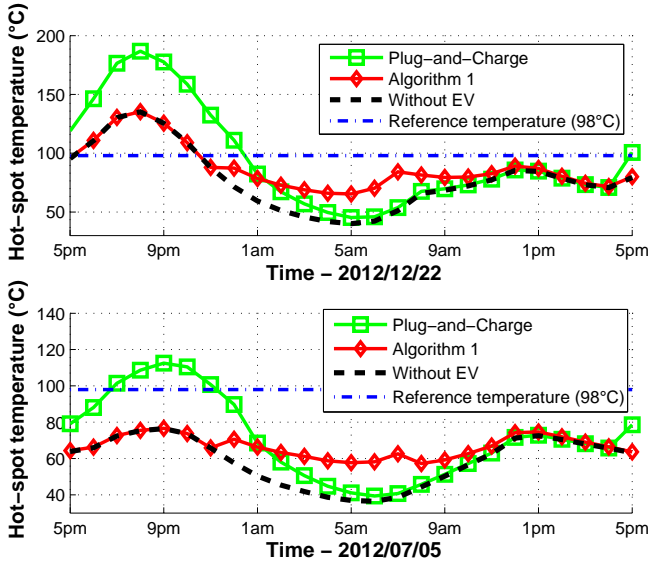}
\caption{Hot-spot temperature against time without EV (for two extreme days in 2012), with Algorithm 1, and with the PaC policy. Scheduling charging start times according to Algorithm 1 instead of plugging-and-charging allows the temperature variations to be minimized and does not induce any excess for the peak. It is also seen that scheduling charging needs properly is much more important during winter time (in France).} 
\label{fig:HSTW}
\vspace{-8mm}
\end{center}
\end{figure}

\textit{HS temperature evolution.} For Scenario (s), Fig. \ref{fig:HSTW} represents the HS temperature against time for the worst day of the year 2012 in France in terms of HS temperature peak (top curve) and the best day (bottom curve) in three scenarios: without EV, with Algorithm 1, and with the PaC policy. The selected worst day (which is a winter day) shows that the HS temperature takes excessive values when the PaC policy is used. The temperature peak is increased by about $80$\degre C w.r.t. to the case without EV (see \cite{watson91}  for physical justifications on the value of this excess). On the other hand, implementing Algorithm 1 does not induce any increase for the peak and roughly tends to minimize HS temperature variations. Since the aging acceleration factor $A_t$ is exponential in the HS temperature, the PaC policy has a dramatic effect in terms of transformer lifetime. Note that this decrease might be made much slower by choosing a transformer with a larger admissible active power (e.g., $120$ kW instead of $90$ kW) but here, our goal is to see the impact of EVs on residential transformers which are already deployed and have therefore been sized to operate without EV.

\begin{figure}
\begin{center}
\includegraphics[scale=0.66]{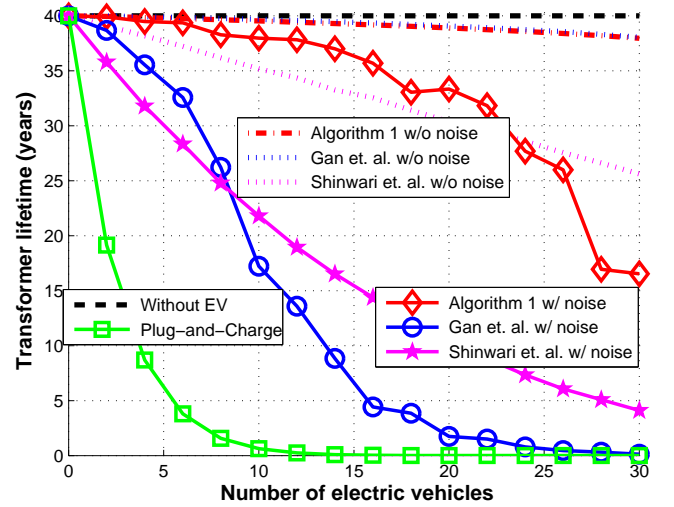}
\caption{Transformer lifetime against the number of EVs ($I$). The plug-and-charge policy is seen to be non-acceptable. Algorithm 1 is seen to perform as well as existing valley-filling type solutions under perfect forecasting of the non-EV demand and outperforms these solutions under imperfect forecasting.}
\label{fig:lifeOfI}
\vspace{-8mm}
\end{center}
\end{figure}

\begin{figure}
\begin{center}
\includegraphics[scale=0.65]{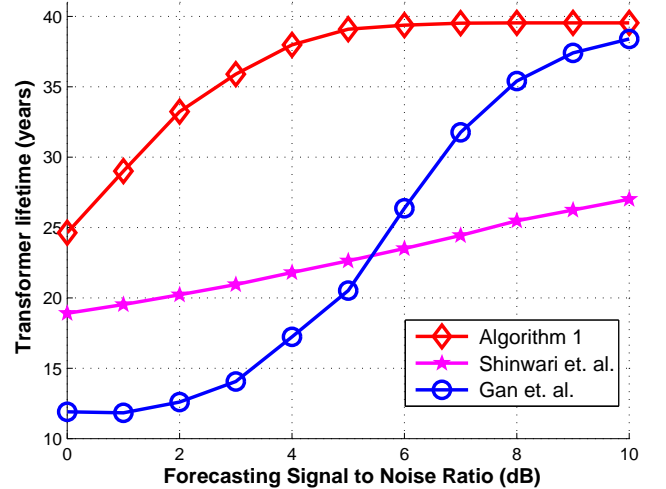}
\caption{Transformer lifetime (years) against the forecasting signal-to-noise ratio (dB) for $I=10$ electric vehicles; $\alpha=1$ and mobility data of Scenario (t), namely those from \cite{entd}, are assumed.}
\label{fig:lifeOfFsnr}
\vspace{-6mm}
\end{center}
\end{figure}

\textit{Performance comparison analysis and influence of the forecasting noise.} Here, we set $\alpha$ to one. The transformer is assumed to be chosen to be able to operate for $40$ years without EV. For Scenario (t), Fig. \ref{fig:lifeOfI} represents the transformer lifetime against the number of EVs  when the exogenous demand forecast is perfect (see the three curves in dotted lines) and when $\mathrm{FSNR}=4$ dB (see the three curves in solid lines). The top curve is an horizontal line which corresponds to the case without EV; the bottom curve corresponds to the PaC policy. The three non-trivial charging policies under consideration are the one corresponding to Algorithm 1, that of Gan et al \cite{Gan13} and Shinwari et al \cite{Shinwari2012}. The Gan et al policy corresponds to the convergence point of an iterative algorithm which aims at minimizing a cost which results of two terms: if maximized alone and assuming convergence, the first term would lead to a VF solution; a second term whose role is to stabilize the parallel implementation-based iterative algorithm. The weight put on the latter penalty term ($0.5$) is tuned optimally for Fig. \ref{fig:lifeOfI}; if this weight is not tuned properly, the implementation of Gan et al may lead to significant performance losses \cite{Xi14}. For the Shinwari et al policy, the energy need of EV $i$ is spread by filling the "holes" of the exogenous demand. For each EV, a proportion of the energy needed is allocated to a given time-slot proportionally to $\frac{\delta_t}{\sum_{t=1}^T \delta_t}$ with $\delta_t =- L_t^{\mathrm{exo}} + \max_t L_t^{\mathrm{exo}}$, and the remainder is uniformly allocated. Fig. \ref{fig:lifeOfI} shows that the PaC policy is seen to be non-acceptable, showing the imperious need for advanced charging schemes. It is seen that Algorithm 1, which is based on rectangular charging profiles, performs as well as the continuous power level-based scheme of Gan et al. Both schemes yield a relatively small decrease of transformer lifetime in Scenario (t). This holds under the assumption of perfect forecasting for the non-EV demand profile. However, under the more realistic assumption of imperfect forecasting (FSNR $= 4$ dB), as seen in Fig. \ref{fig:lifeOfI}, transformer lifetime can severely be degraded as the number of EVs increases for the two latter schemes. Algorithm 1, which is based on rectangular charging profiles, is seen to be much more robust against noise on the exogenous demand than VF solutions. This observation is clearly confirmed by Fig.~\ref{fig:lifeOfFsnr} which is also obtained in Scenario (t) and represents the transformer lifetime against FSNR for three charging schemes. Rectangular profiles have the advantage to be less sensitive to amplitude errors than the VF solution since the sole parameter to be tuned is the charging start time. \tcb{We have checked that this message is not mobility data-dependent. Indeed, simulations which are not provided here and exploit the (US) mobility data of \cite{lee12} and \cite{akhavanhejavi14} confirm that rectangular charging profiles are robust against forecasting errors.} This provides a very strong argument in favor of using \tcb{rectangular profiles}. A parallel with the problem of noise robustness for high-order modulations in digital communications can be drawn and an optimal power level might be identified, which is the purpose of the next paragraph.

\begin{figure}
\begin{center}
\includegraphics[scale=0.68]{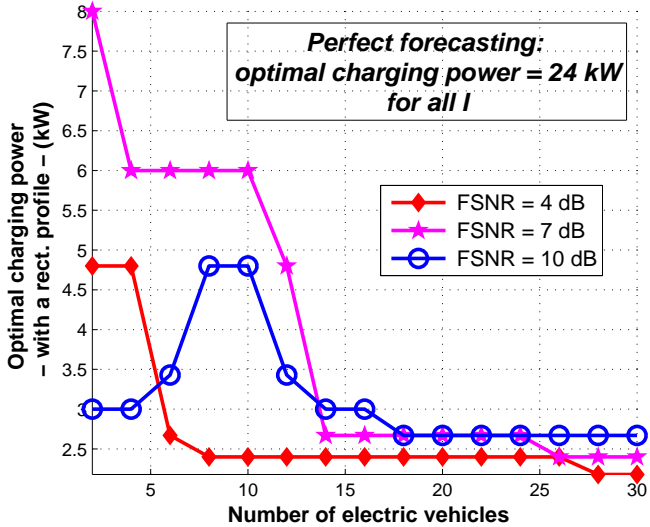}
\caption{In the presence of forecasting noise, charging at the highest charging power possible is suboptimal in terms of transformer lifetime. Rather, a non-trivial optimal charging power level can be determined.}
\label{fig:opt-charging-power}
\vspace{-8mm}
\end{center}
\end{figure}

\textit{Existence of an optimal charging power level for rectangular profiles.}  While it is clear that the larger the charging power the lower the time to charge, having a high charging power can be suboptimal in the presence of forecasting noise. This is what Fig.~\ref{fig:opt-charging-power} shows in Scenario (s). It depicts the optimal charging power in terms of transformer lifetime at the NE obtained with Algorithm 1 as a function of the number of EVs, the minimum (resp. maximum) charging power being set to $2.2$ kW (resp. $24$ kW). At one extreme, in the absence of forecasting noise, the optimal power level is always $24 $ kW, which corresponds to two $30-$min time-slots. At the other extreme, when the noise level is large, the best power level would correspond to charge over the entire period namely at $ \frac{24 \text{kWh}}{15 \text{h}}= 1.6$ kW if it was allowed and we find the minimal value of $2.2$ kW. For typical forecasting noise levels \cite{Kim11}, simulations reveal non-trivial optimal charging power levels, as illustrated by Fig.~\ref{fig:opt-charging-power}.

\begin{figure}[h!]
 \vspace{-4mm}
 \begin{center}
  \centering
      \includegraphics[scale=0.62]{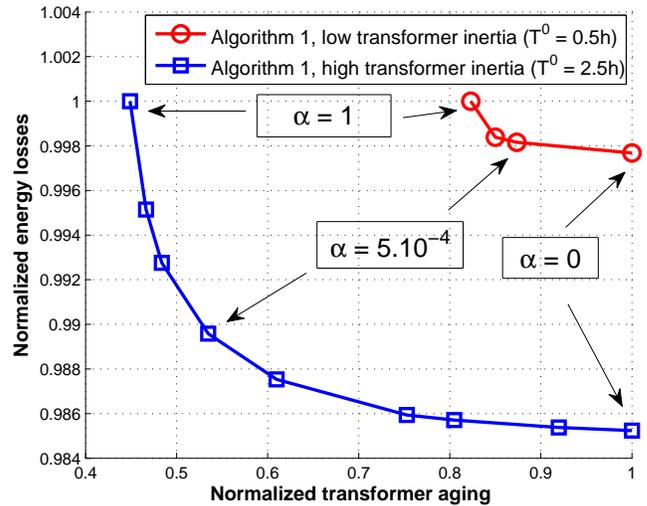}
  \caption{\label{fig:paretoFront} Charging start times have an important impact in terms of transformer aging while they have much less influence on \tcb{energy} losses.} 
  \vspace{-5mm}
  \end{center}
\end{figure}

\textit{Tradeoff between the transformer lifetime and \tcb{energy} losses.} Fig. \ref{fig:paretoFront} represents the Pareto frontier of the feasible cost region for the first day of 2012 when the proposed algorithm is used: the x-axis corresponds to the normalized transformer aging while the y-axis corresponds to normalized \tcb{energy} losses. The curves are obtained by varying $\alpha$ from $0$ to $1$ and considering two different scenarios in terms of thermal inertia: $T^0=0.5$ h (top curve) and $T^0=2.5$ h (bottom curve). The conclusion is that it is preferable to design charging policies which minimize the
transformer aging, i.e., to set $\alpha = 1$. The loss of optimality in terms of \tcb{energy} losses will be rather small by using these policies. It has been observed that changing the charging start time by one hour or two does only affect \tcb{energy} losses in a marginal way; Tab. \ref{tab:jouleOfI}, which provides normalized \tcb{energy} losses for different charging schemes and number of EVs, confirms this. However, because of thermal inertia, changing the start time by one hour or more has a significant impact on aging. This is due to the fact that typical exogenous demand profiles comprise a valley in the night, explaining the results of Fig. \ref{fig:paretoFront}. 
 
\begin{table}
\begin{center}
\begin{tabular}{|c|c|c|c|c|}
\hline
\backslashbox {$J$}{$I$} & $5$ & $10$ & $20$ & $30$ \\
\hline


Plug-and-Charge & $1.14$ & $1.30$ & $1.70$ & $2.18$  \\
\hline
Algorithm 1 & $1.09$ & $1.21$ & $1.50$ & $1.86$  \\
\hline
Gan et. al. \cite{Gan13} & $1.09$ & $1.20$ & $1.49$ & $1.84$  \\
\hline
Shinwari et. al. \cite{Shinwari2012} & $1.10$ & $1.22$ & $1.50$ & $1.85$  \\

\hline
\end{tabular}
\caption{Normalized \tcb{energy} losses $J$ (they are normalized relatively to the case without EV) against the number of electric vehicles ($I$) for four charging schemes.}
\label{tab:jouleOfI}
\vspace{-7mm}
\end{center}
\end{table}

\begin{figure}
\begin{center}
\includegraphics[scale=0.65]{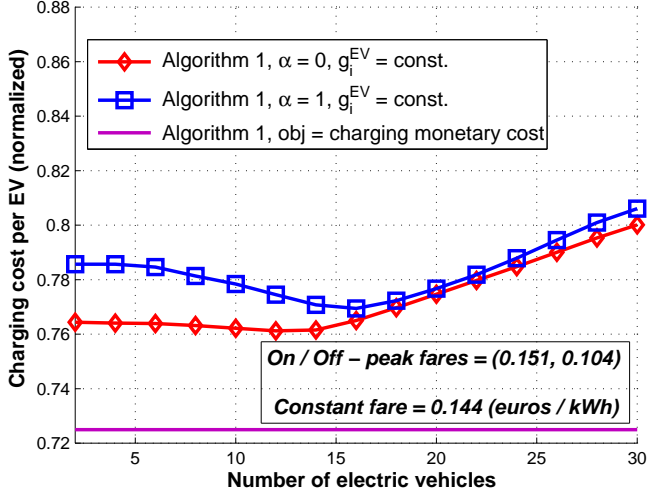}
\caption{EV charging monetary costs at NE (converging point of Algorithm 1) against the number of EVs ($I$) when the individual payoff function \eqref{PersonalCost} corresponds to: \tcb{energy} losses ($\alpha =0$ and $g^{\mathrm{EV}}_{i}=\mathrm{const.}$); transformer aging ($\alpha=1$, $g^{\mathrm{EV}}_{i}=\mathrm{const.}$); the charging monetary cost ($g^{\mathrm{DN}}_{i}=0$ and $g^{\mathrm{EV}}_{i}=\sum_{t=s_i}^{s_i+C_i-1}\pi_{i,t}$). French on/off peak fares (for $\pi_{i,t}$) and Scenario (s) ($C_i=C=16$, $a_i=a=1$, $d_i=d=30$) are assumed.} 
\label{fig:chargingCostHPHC}
\vspace{-8mm}
\end{center}
\end{figure}
\begin{figure}
\begin{center}
\includegraphics[scale=0.65]{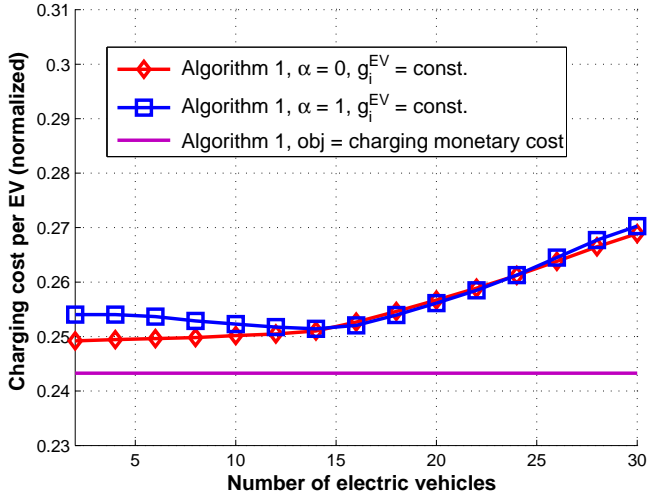}
\caption{This figure correspond to the same setting as Fig.~\ref{fig:chargingCostHPHC} except for the price model which is now the Epex Spot price model (\cite{epex}).} 
\label{fig:chargingCostEpexSpot}
\vspace{-6mm}
\end{center}
\end{figure}

\textit{Tradeoff between transformer aging/\tcb{energy} losses and charging monetary cost.}  The purpose of Fig. \ref{fig:chargingCostHPHC} and Fig. \ref{fig:chargingCostEpexSpot} is to assess what is lost in terms of charging monetary cost when pursuing transformer aging or \tcb{energy} losses minimization. Indeed, these figures depict the charging monetary cost against the number of EVs when: an EV aims at minimizing the charging monetary cost (we force $g_i^{\mathrm{DN}}$ to be zero in \eqref{PersonalCost} and only exploit the function $g_i^{\mathrm{EV}}$); or it aims at minimizing the transformer aging (i.e., when $\alpha=1$ in \eqref{GridCost} and $g_i^{\mathrm{EV}} = \mathrm{const.}$); or it aims at minimizing \tcb{energy} losses  (i.e., when $\alpha=0$ in \eqref{GridCost} and $g_i^{\mathrm{EV}} = \mathrm{const.}$). This is done for two choices of electricity fares: the French on/off peak fare; the market price  (Epex Spot prices in France \cite{epex}). It can be seen that choosing a good charging scheme in terms of aging or \tcb{energy} losses leads to a charging monetary cost which is reasonably close to the one which is obtained by minimizing the monetary cost. The explanation for this is as follows. If the local demand (e.g., at the residential scale) is correlated to the global demand (at the scale of the country), a good tradeoff between technological costs and the monetary cost can be found. Indeed, when looking at the non-EV demand and the electricity price profile, it is seen in  particular that there is an "evening peak" both for the demand and the price. When using other data like the ERCOT (Electricity Reliability Council of Texas, \cite{ercot}) data, similar conclusions can be drawn. Now, if the non-EV demand and prices are not correlated at all, the considered cost function has to combine three terms and the value of the weight $\beta$ (relatively to $\alpha$) will be very influential on the charging profiles obtained with Algorithm 1. Corresponding simulations might then be provided but as our initial goal was to put the emphasis on the distribution network costs and not on the monetary aspect, these simulations are not provided here.

\begin{figure}
\begin{center}
\vspace{-1mm}
\includegraphics[scale=0.65]{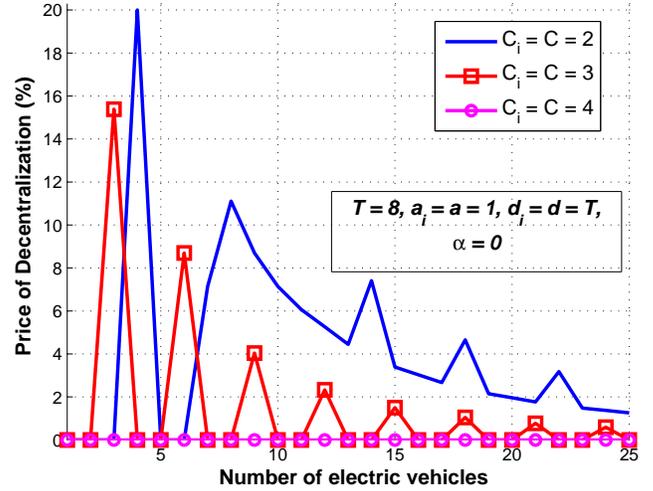}
\caption{Price of Decentralization (PoD) in the simple symmetric case $\mb{L}^{\text{exo}}=0$, $\alpha=0$ (\tcb{energy} losses only), $a_i=a=1$, $d_i=d=T$. \textit{The PoD is small (less than \tcb{$20\%$} for $I=4$ EVs and $C_i=C=2$).}} 
\label{fig:PoD}
\vspace{-1.1cm}
\end{center}
\end{figure}

\textit{Price of Decentralization.} Computing the PoD in general is a hard problem since it involves the maximization of the sum-payoff. Nonetheless, all the numerical results we obtained for various special cases allowed us to confirm what Prop. \ref{PropPoA} suggests namely, that the PoD is close to zero. Concerning the case $\alpha=1$, it is seen from Fig. \ref{fig:lifeOfI} that, under perfect forecasting, the transformer lifetime obtained thanks to Algorithm 1 is greater than $37.75$ years while the no-EV upper bound is at $40$ years. This shows that the PoD is necessarily less than $\frac{40}{37.75}-1 \simeq 6\%$. To elaborate further into this direction, Fig. \ref{fig:PoD} is provided in the case where $\alpha=0$. It represents the PoD against the number of EVs for a simplified setting which allows exhaustive search for the centralized solution to be implemented and which is detailed on the figure itself. Since quantities under use are discrete, there is a combinatorial effect which explains the different peaks. But the general tendency is that the PoD is relatively small and typically decreases with the number of EVs. All these observations fully support the relevance of distributed implementations of charging algorithms.

\section{Concluding remarks}
\label{sec:conclusion}


In this section, we summarize a couple of key messages of this work and then provide possible extensions.

One key feature of this work is to consider rectangular profiles for charging electric vehicles (EVs) and study in details the consequences of this choice. One important message is that rectangular profiles have the advantage to be more robust to errors on the non-EV demand forecast than existing charging schemes which rely on continuous power levels (e.g., valley-filling solutions). As far as transformer lifetime is concerned, it is also seen that in the presence of forecasting noise, there exists an optimal power level at which the EVs should charge.

It is seen that in the context of rectangular charging profiles, it is possible to capture in a simple manner the tradeoff between transformer aging, \tcb{energy} losses, and charging monetary cost. As a rule of thumb for the engineer, it is observed that designing a good charging scheme in terms of aging will imply a good scheme in terms of \tcb{energy} losses whereas the converse does not hold in general (it may hold e.g., for a transformer with a very low thermal inertia). When electricity market prices and the non-EV demand are correlated e.g., by the importance of the evening peak, then a good scheme in terms of aging typically also performs well in terms of charging monetary cost. In more complex scenarios where this correlation is not available, the tradeoff may be analyzed through other simulations which are not provided here. 

The distributed formulation of the charging problem in the considered context is shown to be fully relevant. Indeed, the proposed algorithm is shown to be scalable regarding convergence time, complexity, and required information. Concerning the specific issue of convergence, the ordinality potentiality property is proven to be available for key special cases and therefore guarantees convergence. Otherwise, in more general settings, simulations show that convergence is ensured with overwhelming probability. A possible downside of a distributed formulation of the problem is the potential loss of optimality. However, several strong arguments are provided to show that this issue is quite minor when the problem is formulated as in this work.

A relevant extension is to use a dynamical approach which would take as a system state the transformer hot-spot temperature and charging policies as control policies. Since realistic models for the temperature evolution law are typically non-linear, the considered cost function is typically not quadratic, and charging constraints have to be considered, the underlying control problem seems to be non-trivial. Using a dynamical approach is also relevant to account for the memory effect in electricity  prices. Additionally, the problem of robustness against forecasting errors has to be accounted for. Therefore a quite complete mathematical model would be to account both for memory effects (given by the transfomer and the prices) and uncertainty, which would lead to using the advanced tools of stochastic control \cite{Huang15,Donadee14} and games. 

\appendices

\section{Proof of the potential property of the charging game (Prop. IV.3)}
\label{app:proofPot}

\begin{proof}
We start with the case of assumption (a) i.e., $\forall i \in \mc{I}, \, \mathcal{W}_i(s_i) = \mathcal{W}$. In this case, we have the following sequence of statements:

\begin{enumerate}
\item the game in which the payoff function is the same for all EVs and given by $u_i^1(\bm{s})=u^1(\bm{s})=\displaystyle \sum_{t \in \mathcal{W}} \alpha A_{t}\left(\bm{L}^t\left(\bm{s}\right)\right)+(1-\alpha)J\left(L_{t}\left(\bm{s}\right)\right)$
is an (exact) potential game \cite{Monderer1996} (it is a "team game" more precisely);
\item the game in which the payoff function of EV $i$ is $u_i^2(\bm{s})=g^{\mathrm{EV}}_{i}\left(s_{i}\right)$ is an (exact) potential game because the payoff function of each player only depends on its own strategy;
\item the game in which the payoff function of EV $i$ is $u_i^3(\bm{s})=u_i^1(\bm{s})+u_i^2(\bm{s})$ is an (exact) potential game because the payoff functions of the players are the sum of two payoff functions for which the game is an (exact) potential game;
\item the game in which the payoff function of EV $i$ is $u_i^4(\bm{s})=-f_i(u_i^3(\bm{s}))$ is an ordinal potential game because composing the payoff functions of a potential game by strictly decreasing functions leads to an ordinal potential game\footnote{Furthermore, when $f_i=\mathrm{Id}$, the charging game is an exact potential game.}.
\end{enumerate} 


We now treat case (b). We show that
\begin{align}
\label{eq:PotNoInertia}
\Phi_{(b)}\left(\bm{s}\right)  & =  -\displaystyle \sum_{t \in \mathcal{T}} \displaystyle \sum_{v_{t}=0}^{n_{t}(\bm{s})}  \alpha  A_{t}\left(L_{t}(v_t) \right) + (1-\alpha)J\left(L_{t}(v_t) \right)  \nonumber \\
& -\displaystyle \sum_{i\in \mathcal{I}}g^{\mathrm{EV}}_{i}\left(s_{i}\right)
\end{align}
is an ordinal potential for the charging game. Suppose that player $i$ deviates from $s_{i}$ to $s^{'}_{i}$. Assume, for example, that this deviation is rational and brings him a payoff increase of $ -f_{i}\left(g^{\textrm{DN}}_{i}\left(s_{i}^{\prime},s_{-i}\right)+g^{\textrm{EV}}_{i}\left(s_{i}^{\prime}\right)\right)
+f_{i}\left(g^{\textrm{DN}}_{i}\left(\bm{s}\right)+g^{\textrm{EV}}_{i}\left(s_{i}\right)\right) > 0$

We have the following relations
\begin{align}
& f_{i}\left(g^{\textrm{DN}}_{i}\left(\bm{s}\right)+g^{\textrm{EV}}_{i}\left(s_{i}\right)\right)-f_{i}\left(g^{\textrm{DN}}_{i}\left(s_{i}^{\prime},s_{-i}\right)+g^{\textrm{EV}}_{i}\left(s_{i}^{\prime}\right)\right) > 0 \nonumber \\
\Leftrightarrow & \ g^{\textrm{DN}}_{i}\left(\bm{s}\right)-g^{\textrm{DN}}_{i}\left(s_{i}^{\prime},s_{-i}\right)+g^{\textrm{EV}}_{i}\left(s_{i}\right)-g^{\textrm{EV}}_{i}\left(s_{i}^{\prime}\right) > 0 \label{Equiv1}
\\
\Leftrightarrow & \sum_{t=s_{i}}^{s_{i}+C_{i}-1} g^{\mathrm{DN}}\left(L_{t}\left(\bm{s}\right)\right)-\sum_{t=s_{i}^{\prime}}^{s_{i}^{\prime}+C_{i}-1} g^{\mathrm{DN}}\left(L_{t}\left(\bm{s^{\prime}}\right)\right) \nonumber \\
 & +g^{\textrm{EV}}_{i}\left(s_{i}\right)- g^{\textrm{EV}}_{i}\left(s_{i}^{\prime}\right) > 0 \label{Equiv2} \\
\Leftrightarrow & \, \, \Phi_{(b)}\left(\bm{s}^{\prime}\right)-\Phi_{(b)}\left(\bm{s}\right) >0 \ , \label{Equiv3}
\end{align}
where $g^{\mathrm{DN}}\left(L_{t}\left(\bm{s}\right)\right)=\alpha A_{t}\left(L_{t}\left(\bm{s}\right)\right)
      +(1-\alpha)J\left(L_{t}\left(\bm{s}\right)\right)$ to simplify the notations.
      
The first equivalence (\ref{Equiv1}) comes from the strict monotonicity of $f_{i}$, while (\ref{Equiv2}) and (\ref{Equiv3}) follow by observing that the only terms changing in $\Phi_{(b)}$ because of player $i$'s deviation are $n_{s_{i}},n_{s_{i}+1},...n_{s_{i}+C_{i}-1}$ (decreasing by one), and $n_{s'_{i}},n_{s'_{i}+1},...,n_{s'_{i}+C_{i}-1}$ (increasing by one), which modifies the corresponding sums in (\ref{eq:PotNoInertia}) by, respectively, subtracting a term for $t=s_{i},s_{i}+1,...s_{i}+C_{i}-1$, and adding a term for $t=s'_{i},s'_{i}+1,...s'_{i}+C_{i}-1$. It therefore gives the equivalence in the Ordinal Potential Game definition (Def. \ref{def:OPG}) and concludes the proof.
\end{proof}
\vspace{-3mm}

\section{Proof of the PoD in a special case}
\label{app:proofPoD}

\begin{proof}
We consider scenario (b) for the proof. The proof for $\alpha=0$ is obtained directly from this one observing that setting $\alpha=0$ leads to a DN cost without inertia. Again, let $g^{\mathrm{DN}}\left(L_{t}\left(\bm{s}\right)\right)=\alpha A_{t}\left(L_{t}\left(\bm{s}\right)\right)
      +(1-\alpha)J\left(L_{t}\left(\bm{s}\right)\right)$ to simplify the notations. By following the steps in the proof of App. \ref{app:proofPot}, it can be shown that the nonatomic version of the charging game is also an OPG, with a potential function candidate
      
\begin{equation}
\label{eq:PotNonatomic}
\Phi^{\textrm{NA}}_{(b)}\left(\bm{x}\right)  =  -\displaystyle \sum_{t \in \mathcal{T}} \displaystyle \int_{v_{t}=0}^{x_{t}} g^{\mathrm{DN}}\left(L_t(v_t)\right)\textrm{d}v_t
\end{equation} 
where $L_t(v_t)=L_t^{\text{exo}} + p v_t \, .$ Note that there is no individual term related to $g^{\textrm{EV}}_{i}$ in (\ref{eq:PotNonatomic}) given the assumption made in Prop. \ref{PropPoA}.

Observe now that without the constraints $x_1 \leq \cdots \leq x_C$ and $x_{T-C+1} \geq \cdots \geq x_T$, which ensure that no EV starts charging before $t=1$ or finishes after $t=T$, the minimization problems of both $-\Phi^{\textrm{NA}}_{(b)}$ and $-w$ on the set $\tilde{\mr{X}}=\lbrace \bm{x} \in \left[0,1 \right]^T, \ \sum_{t=1}^T x_t=C \rbrace$ are standard valley-filling problems\footnote{$\min_{\bm{x}} \sum_t f(L_t^{\text{exo}}+px_t) \text{ subject to } \bm{x} \geq 0, \ \sum_t x_t=C$, with $f$ the cost function.} with a strictly convex cost function, here $x \longmapsto \int_{v=0}^{x} g^{\mathrm{DN}}\left(L_t(v)\right)\textrm{d}v$ in the case of $\Phi^{\textrm{NA}}_{(b)}$ ($g^{\mathrm{DN}}$ is strictly increasing) and $x \longmapsto g^{\textrm{DN}}(L_t(x))$ in the case of $w$.
It is known that the solution of this problem is unique and independent of the strictly convex cost function (see e.g., Thm.2 in \cite{Gan13}). Let $\tilde{\bm{x}}^*$ denote this solution. $\tilde{\bm{x}}^*$ has the following "valley-filling" structure:

\vspace{-2mm}
\begin{equation}
\label{eq:VF}
\begin{cases}
\forall t, \ \tilde{x}^*_t>0 \Rightarrow L_t^{\text{exo}}+p\tilde{x}_t^*=L^* \\
\forall t, \ \tilde{x}^*_t=0 \Rightarrow L_t^{\text{exo}} \geq L^* \\
\end{cases}
\ ,
\end{equation}
where $L^*$ denotes the "valley level" of the solution. This implies here that the solutions of the minimization problems of both $\Phi^{\textrm{NA}}_{(b)}$ and $w$ on $\tilde{\mr{X}}$ coincide. 

It remains to show that under the assumptions $L_1^{\text{exo}} \geq \cdots \geq L_C^{\text{exo}}$ and $L_{T-C+1}^{\text{exo}} \leq \cdots \leq L_T^{\text{exo}}$, $\tilde{\bm{x}}^*$ is still the solution of both problems on $\mr{X}$ (instead of $\tilde{\mr{X}}$), which is equivalent to $\tilde{\bm{x}}^* \in \mr{X}$. In fact, it is easy to see that $\tilde{x}^*_1 \leq \cdots \leq \tilde{x}^*_C$ and $\tilde{x}^*_{T-C+1} \leq \cdots \leq \tilde{x}^*_T$ necessarily hold. Suppose indeed that the converse is true e.g., $\tilde{x}^*_1 > \tilde{x}^*_2$. This implies $L_1^{\text{exo}}+p\tilde{x}_1^*>L_2^{\text{exo}}+p\tilde{x}_2^*$ which in turn implies $\tilde{x}^*_1=0$ because of (\ref{eq:VF}), which is contradictory. Then, $\tilde{\bm{x}}^*$ is the (unique) solution of both optimization problems on $\mr{X}$. Thus, the (unique) NE profile is maximizing $w$ on $\mr{X}$ and $\mr{PoD}=0$.
\end{proof}

\ifCLASSOPTIONcaptionsoff
  \newpage
\fi

\bibliographystyle{IEEEtran_NoURL_V2}
\bibliography{LittTransacSG_V11}

\end{document}